\documentclass{AIMS}
 \usepackage{amsmath}
 \usepackage{graphics} 
 \usepackage{epsfig} 
\usepackage[logonly]{trace}

\usepackage{amsfonts, amsmath, amssymb}
\usepackage{graphicx}
\usepackage{pstricks}
\usepackage{psfig}
\usepackage{multirow}
\usepackage{fancyhdr}
\usepackage{vmargin, fancybox}

\usepackage{amsmath,amssymb,amsthm}
\usepackage{t1enc}
  \usepackage{hyperref} 

\def\currentvolume{X}
 \def\currentissue{X}
  \def\currentyear{200X}

   \def\currentmonth{XX}
    \def\ppages{X--XX}

\newtheorem{theorem}{Theorem}[section]

\theoremstyle{definition}
\newtheorem{definition}[theorem]{Definition}

\newtheorem{example}[theorem]{Example}

\title[Towards new schemes: A Lie-group approach
\\of the CBKDV and its derived
equations]
      {Towards new schemes: A Lie-group approach\\ of the CBKDV and its derived equations
      }

\author[Emma Hoarau, Claire David]{}


\email{david@lmm.jussieu.fr}

\begin{document}

\maketitle

\centerline{\scshape Emma Hoarau $\dag$, Claire David $\dag\dag$}
\medskip
{\footnotesize
 \centerline{$\dag$ ONERA, D\'epartement de Simulation Num\'erique et A\'eroacoustique (DSNA),}
 \centerline{ BP 72, 29 avenue de la division Leclerc}
 \centerline{92322 Ch\text{\^ a}tillon Cedex, France}
 \centerline{$\dag\dag$ Universit\'e Pierre et Marie Curie-Paris 6}
  \centerline{Institut Jean Le Rond d'Alembert, UMR CNRS 7190}
   \centerline{Bo\^ite courrier $n^0 162$}
   \centerline{4 place Jussieu, 75252 Paris,
cedex 05, France}
} 

\medskip

 \centerline{(Communicated by Aim Sciences)}
 \medskip

\begin{abstract}

The aim of this paper is to propose methods that enable us to build new numerical schemes, which
preserve the Lie symmetries of the original differential equations. To this purpose, the compound
Burgers-Korteweg-de Vries (\textit{CBKDV}) equation is considered. The particular case of the Burgers
equation is taken as a numerical example, and the resulting semi-invariant scheme is exposed.

\end{abstract}

\section{Introduction}
\label{sec:intro}

\noindent Finite difference equations used to approximate the solutions of a differential equation
generally do not respect the symmetries of the original equation, and can lead to inaccurate numerical results. Usually, specific equations are considered, for
which the authors build a scheme preserving the symmetries of the original differential equation.
Yet, it is more interesting to directly consider a class of differential equations, in order to obtain
more general results.\\

\noindent Using the work of Yanenko \cite{Yanenko} and Shokin \cite{Shokin}, who applied the Lie group
theory to finite difference equations by means of the differential approximation, in conjunction with
the approach of Ames et al. \cite{Ames}, we generalize results developed in \cite{HoarauDavid},
\cite{HoarauDavid2}, and expose the invariance condition for a differential approximation of the
\textit{CBKDV} equation. A numerical example, in the particular case of the Burgers equation, is
developed.

\section{Notion of symmetry}

\subsection{Definitions}
\label{Def-Lie-group}

\begin{definition}
A $r$-parameter Lie group is a $r$-dimensional smooth manifold $G_r$, which has the group properties,
such that the group operation of multiplication and inversion are smooth maps.
\end{definition}
\noindent Especially, a Lie group $G_r$ is defined as a group of continuous transformations which act
on an open subset of the Euclidean space $\mathbb{R}^k$ of variables, which change under
the action of $G_r$. We presently concentrate on a local group, the transformations of which are
close to the identity transformation.
\begin{definition}
A $r$-parameter Lie group $G_r$ is a group of point transformations, which acts on $X\times U$, the
space of the independent variables and the dependent ones:
\footnotesize{
\begin{equation}
\displaystyle G_r=\{x_i^{*}=\phi_i(x,u,a);\
u_j^{*}=\varphi_j(x,u,a),\ i=1,\dots,m;\ j=1,\dots,n;a=(a_1,\dots,a_r)\}
\label{eqn:ga}
\end{equation}
} \normalsize
\noindent where $x\in X\subset \mathbb{R}^m$ and $u\in U\subset \mathbb{R}^n$.

\end{definition}

\noindent $G_r$ locally satisfies the group axioms: existence of an identity element, associativity,
inversibility, closure under the binary composition operation. The transformation corresponding to
a zero parameter is the identity transformation.

\noindent Expand the transformations by means of a Taylor series at the zero value of the parameter
$a$:
\footnotesize{
\begin{eqnarray}
x_i^{*}&=\displaystyle{x_i+a_\alpha\frac{\partial \phi_i}{\partial a_\alpha}\Big|_{a=0}+\mathcal{O}(a_\alpha^2),\ \alpha=1,\dots,r}\nonumber \\
u_j^{*}&=\displaystyle{u_j+a_\alpha\frac{\partial
\varphi_j}{\partial a_\alpha}\Big|_{a=0}+\mathcal{O}(a_\alpha^2),\
\alpha=1,\dots,r}
\end{eqnarray}
} \normalsize
\noindent The derivatives of $\phi_i$ and $\varphi_j$ with respect to the parameter $a_\alpha$ are
smooth functions, called \emph{infinitesimals of the group} $G_r$. Denote by $\xi^\alpha_i$ and
$\eta^\alpha_j$ the infinitesimals of $G_r$.

\noindent The point transformation group $G_r$ can be represented by means of the operator
$\mathbf{L_\alpha}$:
\footnotesize{
\begin{equation}
\displaystyle {\mathbf{L_\alpha}=\xi^\alpha_i(x,u)\frac{\partial
}{\partial x_i}+\eta^\alpha_j(x,u) \frac{\partial }{\partial u_j},\
\ \ i=1,\dots,m;\ j=1,\dots,n;\ \alpha=1,\dots,r} \label{eqn:opL}
\end{equation}
} \normalsize
\noindent The operators $\mathbf{L_\alpha}$, $\alpha=1,\dots,r$ are called the \emph{infinitesimal
operators} of $G_r$.

\noindent $\{\mathbf{L_\alpha}$, $\alpha=1,\dots,r\}$ represents the set of tangent vectors to the
manifold $G_r$, when the zero value is assigned to the parameter $a$. The set is a basis of the
Lie-algebra of the infinitesimal operators of $G_r$, the dimension of which is the same as the one
of the Lie group $G_r$.

\noindent The knowledge of the $\mathbf{L_\alpha}$ enables us to determine the point transformations
of the group $G_r$ by solving the equations:
\footnotesize{
\begin{eqnarray}
\displaystyle{\frac{\partial x_i^{*}}{\partial
a_\alpha}}=\xi^\alpha_i(x^*,u^*),\ \ \displaystyle {\frac{\partial
u_j^{*}}{\partial a_\alpha}}=\eta^\alpha_j(x^*,u^*),\ \
i=1,\dots,m;\ j=1,\dots,n;\ \alpha=1,\dots,r
\label{eqn:GroupFromInf1}
\end{eqnarray}
} \normalsize
\noindent in conjunction with the initial conditions:
\footnotesize{
\begin{eqnarray}
\displaystyle{x_i^{*}\big{|}_{a=0}}=x_i;\ \
\displaystyle{u_j^{*}\big{|}_{a=0}}=u_j
\label{eqn:GroupFromInf2}
\end{eqnarray}
} \normalsize

\begin{example}\textbf{The Galilean transformation group}

\noindent The Galilean transformations correspond to time dependent translations of a reference
frame :
$$G = \left\{\mathcal{T}:(x,t,u)\mapsto(x^*,t^*,u^*)=(x+\epsilon t,t,u+\epsilon)\right\}$$
\noindent where $\epsilon$ is the translation constant velocity.

\noindent $\mathcal{T}$ and its inverse function are continuous. The infinitesimals functions of the
group are $(\xi^\epsilon_x,\xi^\epsilon_t,\eta^\epsilon_u)=(t,0,1)$
\end{example}
\subsection{Symmetry properties of differential equations}

\noindent The notion of symmetry is a tool for generating new solutions of differential equations.
Let us review the main aspects of the application of the Lie group theory to differential equations.

\noindent Consider a system of $l^{th}$-order differential equations:
\footnotesize{
\begin{equation}
\displaystyle \mathcal
F^{\lambda}\big(x,u,u^{(k_1)},u^{(k_1,k_2)},\dots,u^{(k_1\dots
k_l)}\big)=0,\ \ \lambda=1,\dots,q \label{eqn:ED}
\end{equation}
} \normalsize

\noindent Denote by $u^{(k_1\dots k_p)}$ the vector, the components of which are partial derivatives
of order $p$, namely, $u^{(k_1\dots k_p)}_j=\frac{\partial^p u_j}{\partial x_{k_1}\dots\partial
x_{k_p}}$ $j=1,\dots,n$ and $k_1,\dots, k_p \in \{1,\dots,m\}$.

\noindent Denote by $x=(x_1,\dots,x_m)$ the independent variables, $u=(u_1,\dots,u_n)$ the dependent
variables, and $(x_{k_1}\dots x_{k_p})$ a set of elements of the independent variables.

\noindent Equation (\ref{eqn:ED}) is a subset of $X\times
U^{(l)}$, a prolongation of the space $X\times U$ to the space of
the partial derivatives of $u$ with respect to $x$ up to order
$l$. $X\times U^{(l)}$, which is a smooth manifold, is called the
$l$-th order jet space of $X\times U$. In order to take into
account the derivative terms involved in the differential
equation, the action of $G_r$ on $X\times U$ needs to be prolonged
to the space of the derivatives of the dependent variables.

\noindent Denote by $\widetilde G^{(l)}_r$ a $r$-parameter Lie group of point transformations acting
on an open subset $M^{(l)}$ of the $l$-th order jet space $X\times U^{(l)}$ of the independent
variables $x$, dependent variables $u$ and the partial derivatives of $u$ with respect to $x$.

\noindent The $l^{th}$-prolongation operator of $G_r$ is:
\footnotesize{
\begin{eqnarray}
\displaystyle\widetilde{\mathbf{L}}_\alpha^{(l)}=\xi^\alpha_i(x,u)\frac{\partial
}{\partial x_i}+\eta^\alpha_j(x,u) \frac{\partial }{\partial
u_j}+\sigma^{\alpha,(k_1)}_j \frac{\partial }{\partial
{u_j}^{(k_1)}}+\dots+\sigma^{\alpha,(k_1\dots k_l)}_j
\frac{\partial
}{\partial {u_j}^{(k_1\dots k_l)}},\\
\displaystyle \ \ i=1,\dots,m;\ j=1,\dots,n;\
\alpha=1,\dots,r.\nonumber
\end{eqnarray}
} \normalsize

\noindent The infinitesimal functions $\xi^\alpha_i$, $\eta^\alpha_j$, $\sigma^{\alpha,(k_1)}_j$ and
$\sigma^{\alpha,(k_1\dots k_o)}_j$ are given by:
\footnotesize{
\begin{eqnarray}
&&\xi^\alpha_i= \frac{\partial \phi_i}{\partial a_\alpha}\Big|_{a=0},\ \ \ \
\eta^\alpha_j=\frac{\partial\varphi_j}{\partial a_\alpha}\Big|_{a=0},\ \ \ \
\sigma^{\alpha,(k_1)}_j=\displaystyle{\frac{\mathcal{D}\eta^\alpha_j}{\mathcal{D} x_{k_1}}
-\sum_{i=1}^m \frac{\partial u_j}{\partial x_i} \frac{\mathcal{D}\xi^\alpha_i}{\mathcal{D}x_{k_1}}}
\label{eqn:etas}\\
&&\sigma^{\alpha,(k_1\dots k_o)}_j =\displaystyle{\frac{\mathcal{D}
\sigma^{\alpha,(k_1\dots k_{o-1})}_j}{\mathcal{D} x_{k_o}^{\ \ \ \ }}
-\sum_{i=1}^m \frac{\partial^{o} u_j}{\partial x_i \partial x_{k_1}\dots \partial x_{k_{o-1}}}
\frac{\mathcal{D}\xi^\alpha_i}{\mathcal{D}x_{k_o}},\ \ o=2,\dots,l}\nonumber
\end{eqnarray}
} \normalsize

\noindent where:
\footnotesize{$\,\displaystyle\frac{\mathcal{D}}{\mathcal{D}x_k}=\frac{\partial}{\partial x_k}
+\sum_{j=1}^{n}\frac{\partial u_j}{\partial x_k} \frac{\partial }{\partial u_j}$} \normalsize

\begin{theorem}
The system of $l^{th}$-order differential equations is
invariant under the group $\widetilde G^{(l)}_r$ if and only if it satisfies
the following infinitesimal invariance criterion:
\footnotesize{
\begin{equation}
\displaystyle \widetilde{\mathbf{L}}_\alpha^{(l)}\mathcal
F^\lambda \Big|_{\mathcal{F}^{\lambda}=0}=0,\ \ \
\alpha=1,\dots,r;\ \lambda=1,\dots,q \label{eqn:deteq1}
\end{equation}
} \normalsize \label{th:invariance1}
\end{theorem}

\section{Lie group of a differential approximation}

\noindent The symmetry group analysis of the differential
approximation uses the techniques of the Lie group theory applied
to differential equations. The differential approximations involve
step size variable which change under the group action.

\noindent The technique of symmetry analysis is not directly
applied to the finite difference schemes. It is based on the
differential approximation, which describes approximately the
numerical solution behaviour at a reference point of the mesh.
Thus the concept of differential approximation is a local object,
which can not systematically detect a mesh change. Despite the
local behaviour of the symmetry group analysis, the invariant
method based on the differential approximation has enabled one to
establish interesting properties of symmetry group of finite
difference schemes.

\noindent The finite difference scheme, which approximates the differential system (\ref{eqn:ED}),
can be written as:
\footnotesize{
\begin{equation}
\displaystyle{\Lambda^{\lambda}(x,u,h,Tu)=0,\ \ \
\lambda=1,\dots,q} \label{eqn:scheme}
\end{equation}
} \normalsize \noindent where $h=(h_1,h_2,\dots,h_m)$ denotes the space step vector, and
$T=(T_1,T_2,\dots,T_m)$ the shift-operator along the axis of the independent variables, defined by:
\footnotesize{
\begin{equation}
T_i[u](x_1,x_2,\dots,x_{i-1},x_i,x_{i+1},\dots,x_m)=u(x_1,x_2,\dots,x_{i-1},x_i+h_i,x_{i+1},\dots,x_m).
\end{equation}
} \normalsize
\begin{definition}
The differential system: \footnotesize{
\begin{eqnarray}
\displaystyle{\mathcal{P}^{\lambda}\big(x,u,u^{(k_1)},\dots,u^{(k_1\dots k_{l'})}\big)}&=& \displaystyle{\mathcal{F}^{\lambda}\big(x,u,u^{(k_1)},\dots,u^{(k_1\dots k_l)}\big)}\nonumber \\
& & +\displaystyle{\sum_{\beta=1}^s\sum_{i=1}^m (h_i)^{l_\beta} \mathcal{R}^{\lambda}_i(x,u,u^{(k_1)},\dots,u^{(k_1\dots k_{{l'}_{\lambda,i}})})},\nonumber \\
& & \displaystyle{\lambda=1,\dots,q};\ l'= max_{(\lambda,i)}
{l'}_{\lambda,i} \label{eqn:diffapprox}
\end{eqnarray}
} \normalsize is called the $s^{th}$-order differential approximation of the finite difference scheme
(\ref{eqn:scheme}). In the specific case $s=1$, the above system is called the first differential
approximation.
\label{def:diffapprox}
\end{definition}

\noindent The differential system (\ref{eqn:diffapprox}) is obtained from the algebraic system
(\ref{eqn:scheme}) by applying Taylor series expansion to the components of $T u$ about the point
$x=(x_1,\ \dots,\ x_m)$ and truncating the expansion to a given finite order.

\noindent Denote by $G'_r$ a group of transformations acting on an open subset $M'$ of
$X\times U\times H$ the space of the independent variables, the dependent variables and the step size
variables :
\footnotesize{
\begin{equation}
\displaystyle G'_r=\{x_i^{*}=\phi_i(x,u,a);\
u_j^{*}=\varphi_j(x,u,a);h_i^{*}=\psi_i(x,u,h,a),\ i=1,\dots,m;\
j=1,\dots,n\}
\end{equation}
} \normalsize
\noindent by $\mathbf{L_\alpha}'$ the basis infinitesimal operator of $G'_r$:
\footnotesize{
\begin{equation}
\displaystyle{\mathbf{L_\alpha}'= \xi^\alpha_i(x,u)\frac{\partial
}{\partial x_i}+\eta^\alpha_j(x,u) \frac{\partial }{\partial
u_j}+\zeta^\alpha_i(x,u,h)\frac{\partial}{\partial h_i},\ \ \
\alpha=1,\dots,r}
\end{equation}
} \normalsize where \footnotesize{$\zeta^\alpha_i=\frac{\partial
\psi_i}{\partial a_\alpha}\Big{|}_{a=0},\ \
\alpha=1,\dots,r$}\normalsize

\noindent and by $\widetilde{G}^{(l')}_r$ a group of transformation acting on an open subset
$M^{(l')}$ of the space of the independent variables, the dependent variables and the step size
variables and the partial derivatives involved in the differential system \label{eqn:diffapprox}.
The ${l'}^{th}$-prolongation operator of $G'_r$, $\widetilde{\mathbf{L}}_\alpha^{(l')}$ can be
written as:
\footnotesize{
\begin{equation}
\displaystyle
\widetilde{\mathbf{L}}_\alpha^{(l')}=\mathbf{L_\alpha}'+\sum_{j=1}^n\sum_{p=1}^{l'}\sigma_j^{\alpha,{(k_1\dots
k_p)}}\frac{\partial }{\partial u_j^{(k_1\dots k_p)}}
\end{equation}
} \normalsize
\begin{theorem}
\noindent The differential approximation (\ref{eqn:diffapprox}) is
invariant under the group $\widetilde G^{(l')}_r$ if and only if
\footnotesize{
\begin{eqnarray}
\displaystyle\widetilde{\mathbf{L}}_\alpha^{(l')}\mathcal{P}^{\lambda}\big((x,u,u^{(k_1)},\dots,u^{(k_1\dots
k_{l'})}\big)\Big|_{\mathcal{P}^{\lambda}=0}&=&0,\ \ \
\alpha=1,\dots,r;\ \lambda=1,\dots,q\\
\displaystyle \text{or}\, \,
\Big[\widetilde{\mathbf{L}}_\alpha^{(l)}\mathcal{F}^{\lambda}+\widetilde{\mathbf{L}}_\alpha^{(l')}\Big(\sum_{\beta=1}^s\sum_{i=1}^m
(h_i)^{l_\beta}
\mathcal{R}^{\lambda}_i\Big)\Big]\Big|_{\mathcal{P}^{\lambda}=0}&=&0,\
\ \ \alpha=1,\dots,r;\ \lambda=1,\dots,q \label{eqn:deteq2}
\end{eqnarray}
} \normalsize \label{th:invariance2}
\end{theorem}

\section{Determination of the infinitesimal functions}

\noindent The calculation of Lie groups of differential equations with pencil and paper is tedious and
may induce errors. The size of related equations increases with the number of the symmetry variables,
and the order of the differential equations. A large amount of packages have been created using software
programs with symbolic manipulations, such as \textsf{Mathematica}, \textsf{MACSYMA}, \textsf{Maple},
\textsf{REDUCE}, \textsf{AXIOM}, \textsf{MuPAD}. Schwarz \cite{Schwarz} wrote algorithms for
\textsf{REDUCE} and \textsf{AXIOM} computer algebra systems, Vu and Carminati \cite{VuCarminati} worked
on DESOLVE, a \textsf{Maple} program, Herod \cite{Herod} and Baumann \cite{Baumann} developed
\textsf{Mathematica} programs.

\noindent The authors have implemented a \textsf{Mathematica} package \cite{HoarauDavid}, for the
determination of the Lie group of the differential approximation of one dimensional model equations.

\noindent {Theorems \ref{th:invariance1}} and
{\ref{th:invariance2}}, respectively give the algorithmic
procedure for the determination of Lie group of any differential
system and differential approximation.

\noindent The theorem infinitesimal invariance criteria involved the independent variables $x$, the
dependent ones $u$, products of the partial derivatives of $u$ with respect to $x$, the unknown
infinitesimal functions $\xi^\alpha_i$, and $\eta^\alpha_j$, $i=1,\dots,m$; $j=1,\dots,n$ and their
partial derivatives with respect to $x$ and $u$.

\noindent The partial derivatives of the infinitesimal functions with respect to $x$ and $u$ come from
the coefficients of prolonged Lie algebra vector field (\ref{eqn:etas}).

\noindent Equations (\ref{eqn:deteq1}) and (\ref{eqn:deteq2}) are
simplified by means of the conditions (\ref{eqn:ED}) and
(\ref{eqn:diffapprox}). This simplification manipulation
eliminates some derivatives of $u$. These equations are then
solved algebrically with respect to the partial derivatives of the
dependent variables, handled as independent variables. Denote by
$w$ the vector, the components of which are these variables. Since
the whole equation holds for all the $w$ components, each
coefficient in front of the products of the $w$ components has to
be zero. This leads to a linear overdetermined system of partial
differential equation, with respect to the infinitesimal
functions, called the \textit{determining equations} of the Lie
group of the related differential system .

\noindent The solution of the overdetermined system can be found
either by using elementary methods of the theory of linear partial
differential equations or by using a polynomial form for the
infinitesimals. The last technique provides a linear system of
algebraic equations with respect to the polynomial coefficients.
But it can not find infinitesimals with transcendal functions.

\noindent The resolution of the determining equations yields explicitly the expression of $\xi^{\alpha}_i$,
$\eta^{\alpha}_j$, $\alpha=1,\dots,r,\ i=1,\dots,m,\ j=1,\dots,n$. Then relations (\ref{eqn:GroupFromInf1})
and (\ref{eqn:GroupFromInf2}) provide the calculation of group transformations from the infinitesimal
expression.

\section{Symmetries of the CBKDV equation}

\noindent The invariance condition of {theorem \ref{th:invariance1}}, which enables us to obtain the
expression of the infinitesimal operators, can be written as: \footnotesize{
\begin{equation}
\displaystyle{{\tilde{\mathbf{L}}}_\alpha^{(2)} \mathcal{F}
\Big{|}_{u_t + \alpha u u_x + \beta u^2 u_x + \mu u_{xx} - s
u_{xxx}=0}=0,\, \, \, \alpha=1,\dots,r;} \label{eqn:condinvBurg}
\end{equation}
} \normalsize

\noindent Equation (\ref{eqn:condinvBurg}) provides an overdetermined system of linear partial
differential equations for the infinitesimal functions.

\noindent Solving these equations, by means of a symbolic calculus tool (\textsf{Mathematica}),
provides the expression of the infinitesimal functions of a two-parameter group:
\footnotesize{
\begin{eqnarray}
\displaystyle{\xi^\alpha_1}=a_0\,\,\, ,  \,\,\,
\displaystyle{\xi^\alpha_2=b_0}\,\,\, ,  \,\,\,
\displaystyle{\eta^\alpha=0}\nonumber
\end{eqnarray}
} \normalsize

\noindent where $a_0$, $b_0$ are constants.

\noindent The two-dimensional Lie algebra of the group $G$ is generated by the following operators:
\footnotesize{
\begin{eqnarray}
& &\displaystyle{\mathbf{L}_1=\frac{\partial}{\partial x}},\ \ \
 \displaystyle{\mathbf{L}_2=\frac{\partial}{\partial t}}
\nonumber
\end{eqnarray}
} \normalsize

\noindent which respectively correspond to:
\begin{itemize}
\item the space translation
\footnotesize{$:\,(x,t,u,\nu) \longmapsto
(x+x_0,t,u,\nu)$}\normalsize ;
\item the time translation
\footnotesize{$:\,(x,t,u,\nu) \longmapsto
(x,t+t_0,u,\nu)$}\normalsize ;
\end{itemize}

\normalsize \noindent $x_0,\ t_0$ are constants.\\

\noindent It is important to note that, in the case of particular cases of the CBKDV equation, when,
for instance, some of the coefficients are equal to zero, more symmetries can be obtained.

\section{The specific case of the Burgers equation}

\noindent Yanenko \cite{Yanenko} and Shokin \cite{Shokin} symmetry
analysis is applied to finite difference schemes for solving the
Burgers equation. Thus, the symmetries of the Burgers equation,
which are finally broken by the finite difference discretization,
are determined. The techniques exposed in \cite{Shokin} and
\cite{Yanenko} enable one to construct differential
approximations, which preserve the symmetries of the original
differential system. We call the related finite difference scheme
a \emph{semi-invariant} scheme, in so far as the invariance
condition is weaker than the one of the other invariance methods,
defined as direct invariance methods. Indeed, the approach in
\cite{Shokin} and \cite{Yanenko} does not deal with the invariance
of the algebraic equations, which govern the mesh evolution.

\noindent The differential equation has provided important
characteristics for numerical schemes, in the study of numerical
stability, dissipation and dispersive property. In \cite{Shokin},
\cite{Yanenko} and \cite{Ames}, the differential approximation has
been revealed as a practical and recevable tool for symmetry
analysis of finite difference scheme.

\noindent A comparison is made between the numerical solutions of the Burgers equation for some
standard schemes and the semi-invariant one.

\subsection{Symmetries of the Burgers equation}

\noindent The Burgers equation, which is a particular case of the CBKDV equation, with $\alpha=1$,
$\beta =0$, $s=0$, $\mu=-\nu$, can be written as:
\footnotesize{
\begin{equation}
\displaystyle{\mathcal{F}(x,t,u,\nu,u_x,u_t,u_{xx})=u_t+u\
u_x-\nu\ u_{xx}=0} \label{eqn:burger}
\end{equation}
} \normalsize \noindent where $\nu\geq 0$ is the dynamic viscosity.

\noindent Let us denote by $G$ a group of continuous transformations of the Burgers equation acting
on an open subset $M$ of the space the independent variables $(x,t)$, the dependent variable $u$,
and the viscosity $\nu$. The viscosity is taken as a symmetry variable in order to enable us to take
into account variations of the Reynolds number.

\noindent The six-dimensional Lie algebra of the group $G$ is generated by the following operators:
\footnotesize{
\begin{eqnarray}
& &\displaystyle{\mathbf{L}_1=\frac{\partial}{\partial x}},\ \displaystyle{\mathbf{L}_2=
\frac{\partial}{\partial t}},\ \displaystyle{\mathbf{L}_3=x\frac{\partial}{\partial x}+2 t\frac{\partial}{\partial t}-u\frac{\partial}{\partial u}}\nonumber\\
& &\displaystyle{\mathbf{L}_4=x t\frac{\partial}{\partial
x}+t^2\frac{\partial}{\partial t}+(-u t+x)\frac{\partial}{\partial
u}},\ \displaystyle{\mathbf{L}_5=t\frac{\partial}{\partial
x}+\frac{\partial}{\partial u}},\
\displaystyle{\mathbf{L}_6=-t\frac{\partial}{\partial
t}+u\frac{\partial}{\partial u}+\nu\frac{\partial}{\partial
\nu}}\label{eqn:op-burgers}
\end{eqnarray}
} \normalsize which respectively correspond to:
\begin{itemize}
\item the space translation
\footnotesize{$:\,(x,t,u,\nu) \longmapsto
(x+\epsilon_1,t,u,\nu)$}\normalsize;
\item the time translation
\footnotesize{$:\,(x,t,u,\nu) \longmapsto
(x,t+\epsilon_2,u,\nu)$}\normalsize;
\item the dilatation
\footnotesize{$:\,(x,t,u,\nu) \longmapsto (\epsilon_3
x,\epsilon_3^2 t,\epsilon_3^{-1} u,\nu)$}\normalsize;
\item the projective transformation
\footnotesize{$:\,(x,t,u,\nu) \longmapsto
\Big(\frac{x}{1-\epsilon_4 t},\frac{t}{1-\epsilon_4
t},x\epsilon_4+u(1-\epsilon_4 t),\nu\Big)$}\normalsize;
\item the Galilean transformation
\footnotesize{$:\,(x,t,u,\nu) \longmapsto (x+\epsilon_5\
t,t,u+\epsilon_5,\nu)$}\normalsize;
\item the dilatation
\footnotesize{$:\,(x,t,u,\nu) \longmapsto (x,\epsilon_6^{-1}
t,\epsilon_6 u, \epsilon_6\nu)$}\normalsize.
\end{itemize}
\normalsize {$(\epsilon_i)_{i=1,\dots,6}$} \normalsize are
constants. \normalsize

\noindent As expected, the Burgers equation, as a particular form
of the \textit{CBKDV} equation, admits more symmetries, yieldind a
richer Lie algebra.

\subsection{Symmetries of first differential approximations}
\label{Sym-FDA}

\noindent Denote by $h$ the mesh size, $\tau$ the time step, $N_x$ the number of mesh points, $N_t$
the number of time steps, and $u^n_i,\ i \in \{0,\dots,N_t\},\ n \in \{0,\dots,N_x\}$ the discrete
approximation of $u(i h,n \tau)$.

\noindent In order to shorten the size of the finite difference scheme expressions, we use the
following notations introduced by Hildebrand in \cite{Hildebrand}:
\footnotesize{
\begin{eqnarray*}
\displaystyle \delta (u^n_i)=\frac{u^n_{i+\frac{1}{2}}-u^n_{i-\frac{1}{2}}}{h},& & \mu (u^n_i)=\frac{u^n_{i+\frac{1}{2}}+u^n_{i-\frac{1}{2}}}{2}\\
\displaystyle \delta^{+} (u^n_i)=\frac{u^n_{i+1}-u^n_{i}}{h},& &
\delta^{-} (u^n_i)=\frac{u^n_{i}-u^n_{i-1}}{h},\,\,\,\,\,\,\,\,
E^{\alpha}u^n_{i}=u^n_{i+\alpha}
\end{eqnarray*}
}\normalsize

\noindent The Burgers equation can be discretized by means of:
\begin{itemize}
\item \textbf{the FTCS (forward-time and centered-space) scheme}:
\footnotesize{
$$\displaystyle {\frac{u^{n+1}_i -u^n_i}{\tau}+\frac{\mu\delta}{h}\big(\frac{u^2}{2}\big)^n_{i}-\nu \frac{\delta^2}{h^2}u^n_{i}=0}$$
} \normalsize
\item \textbf{the Lax-Wendroff scheme}:
\footnotesize{
$$\displaystyle {\frac{u^{n+1}_i -u^n_i}{\tau}+\frac{\mu\delta}{h}\big(\frac{u^2}{2}\big)^n_{i}-\nu \frac{\delta^2}{h^2}u^n_{i}+A^n_i=0}$$
} \normalsize
\noindent where:
\footnotesize{
\begin{eqnarray*}
\displaystyle A^n_i=&-&\frac{\tau}{2 h^2}\Big[E^{\frac{1}{2}}u^n_{i}\ \delta^{+}\big(\frac{u^2}{2}\big)^n_{i}-E^{-\frac{1}{2}}u^n_{i}\
\delta^{-}\big(\frac{u^2}{2}\big)^n_{i}\Big]-\frac{\nu^2\tau}{2}\Big[\frac{\delta^4}{h^4}u^n_{i}\Big]\\
&+&\frac{\nu\tau}{2 h^3}\Big[E^{\frac{1}{2}}u^n_{i}\ \delta^{2}(E^{\frac{1}{2}}u^n_{i})-E^{-\frac{1}{2}}u^n_{i}\ \delta^{2}(E^{-\frac{1}{2}}u^n_{i})\Big]
+\frac{\nu\tau}{2}\Big[\frac{\mu \delta^3}{h^3}\big(\frac{u^2}{2}\big)^n_{i}\Big]
\end{eqnarray*}
} \normalsize
\item \textbf{the Crank-Nicolson scheme}:
\footnotesize{
\begin{eqnarray*}
\displaystyle \frac{u^{n+1}_i-u^{n}_i}{\tau}+\frac{\mu\delta}{h}\Big[\big(\frac{u^2}{2}\big)^{n+1}_{i}+\big(\frac{u^2}{2}\big)^n_{i}\Big]-\nu \frac{\delta^2}{h^2}[u^{n+1}_{i}+u^n_{i}]=0
\end{eqnarray*}
} \normalsize
\end{itemize}
\normalsize
\noindent The linear stability properties and the related orders of approximation are:
\begin{itemize}
\item \textbf{the FTCS scheme}: \footnotesize{$S\leq \frac{1}{2}$, $CFL\leq 1$, $CFL^2\leq 2 S$}\normalsize; \footnotesize{$ \mathcal{O}(\tau,h^2)$}\normalsize
\item \textbf{the Lax-Wendroff scheme}: \footnotesize{$S^*\leq \frac{1}{2}$, $CFL\leq 1$}\normalsize; \footnotesize{$\displaystyle \mathcal{O}(\tau^2,h^2)$}\normalsize
\item \textbf{the Crank-Nicolson scheme}: unconditional stability; \footnotesize{$\displaystyle \mathcal{O}(\tau^2,h^2)$}\normalsize
\end{itemize}
\noindent where $CFL=\frac{a\tau}{h}$, $S=\frac{\nu\tau}{h^2}$ and $S^*=\big(\nu+\frac{a h CFL}{2}\big)\frac{\tau}{h^2}$.

\noindent Consider ${u_{i}}^{n}$ as a function of the time step $\tau$, and of the mesh size $h$, expand
it at a given order by means of its Taylor series, and neglect the $o(\tau^{\alpha})$ and $o({h}^{\beta})$
terms, where $\alpha$ and $\beta$ depend on the order of the schemes. This yields the differential
representation of the finite difference equation.

\noindent The following differential representations are obtained:
\begin{itemize}
\item\textbf{for the FTCS scheme:}
\footnotesize{
\begin{eqnarray*}
\displaystyle u_t+\frac{1}{2}(u^2)_x-\nu\
u_{xx}+\frac{\tau}{2}g_2+\frac{h^2}{12}(u^2)_{xxx}-\frac{\nu
h^2}{12}u_{xxxx}=0
\end{eqnarray*}
} \normalsize
\item \textbf{for the Lax-Wendroff scheme:}
\footnotesize{
\begin{eqnarray*}
\displaystyle u_t+\frac{1}{2}(u^2)_x-\nu\
u_{xx}+\frac{\tau^2}{6}g_3+\frac{h^2}{12}(u^2)_{xxx}-\frac{\nu
h^2}{12}u_{xxxx}=0
\end{eqnarray*}
} \normalsize
\item \textbf{for the Crank-Nicolson scheme:}
\footnotesize{
\begin{eqnarray*}
u_t+\frac{1}{2}(u^2)_x-\nu
u_{xx}+\tau^2\big(\frac{g_3}{6}+\frac{1}{4}(g^2_1+ug_2)_x-\frac{\nu}{4}(g_2)_{xx}\big)+\frac{h^2}{12}(u^2)_{xxx}-\frac{\nu h^2}{12}u_{xxxx}=0
\end{eqnarray*}
} \normalsize
\end{itemize}
\noindent where \footnotesize{$g_1=-\big(\frac{u^2}{2}\big)_x+\nu u_{xx}$,
$g_2=\big(-g_1 u\big)_x+\nu\big(g_1\big)_{xx}$, $g_3=\big(-g_2 u
-g^2_1\big)_x+\nu \big(g_2\big)_{xx}$ } \normalsize \vspace{0.4cm}

\noindent Denote by $G'$ the group of transformations of a first differential approximation acting on
an open subset $M'$ of the space of the independent variables $(x,t)$ and the dependent variable $u$,
the step size variables $(h,\tau)$ and the viscosity $\nu$.

\noindent The ${l'}^{th}$-prolongation of $G'$ can be written as:
\footnotesize{
\begin{eqnarray}
\displaystyle{\widetilde{\mathbf{L}}_\alpha'^{(l')}}&=&\displaystyle{\xi^\alpha_1\frac{\partial}{\partial
x}+\xi^\alpha_2\frac{\partial}{\partial t}+\eta^\alpha
\frac{\partial}{\partial
u}+\sum_{p=1}^{l'}\sigma_j^{\alpha,(k_1\dots k_p)}\frac{\partial
}{\partial u_j^{(k_1\dots
k_p)}}+\zeta^\alpha_1\frac{\partial}{\partial
h}+\zeta^\alpha_2\frac{\partial}{\partial
\tau}+\theta^\alpha\frac{\partial}{\partial \nu}}
\end{eqnarray}
}\normalsize \noindent where $l'$ has been defined in
\textbf{definition {\ref{def:diffapprox}}}.

\noindent \textbf{Theorem \ref{th:invariance2}} enables us to
obtain the necessary and sufficient condition of invariance of the
first differential approximation $\mathcal{P}$: \footnotesize{
\begin{equation}
\displaystyle{\widetilde{\mathbf{L}}_\alpha'^{(l')}
\mathcal{P}\Big{|}_{\mathcal{P}=0}=0}
\label{eqn:condinvdiffapprox}
\end{equation}
} \normalsize

\noindent \textbf{Theorem \ref{th:invariance2}} is applied to the differential representations
of the above schemes.

\noindent The resolution of the determining equations of each first differential approximation yields the
$4$-parameter group (see \cite{HoarauDavid}):
\footnotesize{
\begin{eqnarray}
\displaystyle{\xi^\alpha_1=a+b\ x},&\ \ \displaystyle{\xi^\alpha_2=c+(2 b-d)\  t},&\ \ \displaystyle{\eta^\alpha=(-b+d)\ u}\\
\displaystyle{\zeta^\alpha_1=b\ h},&\ \
\displaystyle{\zeta^\alpha_2=(2 b-d)\ \tau},&\ \
\displaystyle{\theta^\alpha=e \nu}\nonumber
\end{eqnarray}
} \normalsize

\noindent The $4$-dimensional Lie algebra of $G'$ is generated by:
\footnotesize{
\begin{eqnarray}
& &\displaystyle{\mathbf{L}_1=\frac{\partial}{\partial x}},\ \
\displaystyle{\mathbf{L}_2=\frac{\partial}{\partial t}}\nonumber ,\ \ \displaystyle{\mathbf{L'}_3=x\frac{\partial}{\partial x}+2 t\frac{\partial}{\partial t}-u\frac{\partial}{\partial u}+h\frac{\partial}{\partial h}+2 \tau\frac{\partial}{\partial \tau}}\\
& & \displaystyle{\mathbf{L'}_4=-t\frac{\partial}{\partial
t}+u\frac{\partial}{\partial u}-\tau \frac{\partial}{\partial
\tau}+\nu \frac{\partial}{\partial \nu}}
\end{eqnarray}
} \normalsize

\noindent These operators are respectively related to:
\begin{itemize}
\item the space translation
\footnotesize{$:\,(x,t,u,h,\tau,\nu) \longmapsto
(x+\epsilon_1,t,u,h,\tau,\nu)$}\normalsize ;
\item the time translation
\footnotesize{$:\,(x,t,u,h,\tau,\nu) \longmapsto
(x,t+\epsilon_2,u,h,\tau,\nu)$}\normalsize ;
\item the dilatation
\footnotesize{$:\,(x,t,u,h,\tau,\nu) \longmapsto (\epsilon_3
x,\epsilon_3^2 t,\epsilon_3^{-1} u,\epsilon_3 h,\epsilon_3^2
\tau,\nu)$}\normalsize ;
\item the dilatation
\footnotesize{$:\,(x,t,u,h,\tau,\nu) \longmapsto
(x,\epsilon_4^{-1} t,\epsilon_4 u,h,\epsilon_4^{-1}\tau,\epsilon_4
\nu)$}\normalsize ;
\end{itemize}
\noindent where $(\epsilon_i)_{i=1,\dots,4}$ are constants.\\
\noindent The above finite difference equations are preserved by
the space translation, the time translation and both dilatations.

\noindent Approximating the Burgers equation by the above finite
difference equations results in the loss of the projective and
Galilean transformations.

\subsection{A semi-invariant scheme}

\noindent Two analogous methods provide a direct symmetry analysis
of finite difference schemes and can lead to the definition of
adapted evolutionary meshes, whose geometrical structure is
preserved by the entire group. The first direct method has been
introduced by Dorodnitsyn \cite{Dorodnitsyn-theory} and is based
on Lie algebra techniques, using the infinitesimal operators. The
second method has been introduced by Olver \cite{Olverinv} and is
based on the theory of the Cartan moving frame. The method
proposed by Yanenko \cite{Yanenko} and Shokin \cite{Shokin}
consists of a symmetry study of the differential approximation.
Although the last method is not fully exact, the numerical results
in \cite{Shokin} and \cite{Ames} has proved its effectiveness.

\noindent The scheme proposed below is associated to an uniform
orthogonal mesh.

\noindent We propose to approximate the Burgers equation by the
finite difference scheme: \footnotesize{
\begin{eqnarray}
\displaystyle\frac{u^{n+1}_i-u^n_i}{\tau}+\frac{1}{h}\big(\mu\delta-\frac{\mu\delta^3}{6}\big)\big(\frac{u^2}{2}\big)^n_{i}-\nu\frac{1}{h^2}\big(\delta^2-\frac{\delta^4}{12}\big)(u^n_i)-h\Big(\Omega^n_{i+\frac{1}{2}}
\delta^{+}-\Omega^n_{i-\frac{1}{2}} \delta^{-}\Big)u^n_i=0
\label{eqn:inv}
\end{eqnarray}
} \normalsize
\noindent where $\Omega^n_i=\Omega(x_i,t_n,u^n_i)$
is defined next so that the related differential representation is
preserved by the symmetries of the Burgers equation. \noindent The
scheme has second-order accuracy in space and first-order accuracy
in time. The derivatives $(u^2)_x$ and $u_{xx}$ are approximated
by fourth order accuracy difference expressions: \footnotesize{
\begin{eqnarray}
\displaystyle\big(\frac{\mu\delta}{h}-\frac{\mu\delta^3}{6
h}\big)(u^n_i)=\big(u_x-\frac{h^4}{30}u_{5x}\big)^n_i+\mathcal{O}(h^6),\
\displaystyle\big(\frac{\delta^2}{h^2}-\frac{\delta^4}{12
h^2}\big)(u^n_i)=\big(u_{xx}-\frac{h^4}{90}u_{6x}\big)^n_i+\mathcal{O}(h^6)
\label{eqn:4-order-acc}
\end{eqnarray}
} \normalsize \noindent The truncation error of the difference
scheme (\ref{eqn:inv}) can be written as: \footnotesize{
\begin{eqnarray}
\displaystyle \epsilon&=&\frac{\tau}{2}u_{tt}-{h^2}\Big(\Omega
u_x\Big)_x+\mathcal{O}(\tau^2)+\mathcal{O}(h^4)\nonumber
\end{eqnarray}
} \normalsize $u_{tt}$ is replaced by an expression involving
partial derivatives with respect to $x$, by using the Burgers
equation. Replacing the obtained expression in the truncation
error leads to: \footnotesize{
\begin{eqnarray*}
\displaystyle \epsilon&=&\Big(C u_x\Big)_x-\frac{\nu\tau}{2}\Big(u
u_{xx}\Big)_x-\frac{\nu\tau}{2}\big(\frac{u^2}{2}\big)_{xxx}+\frac{\nu^2\tau}{2}u_{xxxx}+\mathcal{O}(\tau^2)+\mathcal{O}(h^4)
\end{eqnarray*}
} \normalsize
\noindent where $C=\frac{\tau}{2}u^2-h^2\Omega$.\\
\noindent It is convenient for the calculation of $C$ that the
truncation error is reduced to: \footnotesize{
\begin{eqnarray*}
\displaystyle \epsilon&=&\Big(C
u_x\Big)_x+\mathcal{O}(\tau^2)+\mathcal{O}(h^4)
\end{eqnarray*}
} \normalsize \noindent The related finite difference scheme is
the following first order accuracy in time and second order
accuracy in space: \footnotesize{
\begin{eqnarray}
\displaystyle
\frac{u^{n+1}_i-u^n_i}{\tau}+\frac{1}{h}\big(\mu\delta-\frac{\mu\delta^3}{6}\big)\big(\frac{u^2}{2}\big)^n_{i}-\nu\frac{1}{h^2}\big(\delta^2-\frac{\delta^4}{12}\big)(u^n_i)-h\Big(\Omega^n_{i+\frac{1}{2}}
\delta^{+} -\Omega^n_{i-\frac{1}{2}} \delta^{-}\Big)u^n_i\nonumber\\
\displaystyle
+\frac{\nu\tau}{2}\Big(u^n_{i+\frac{1}{2}}\frac{\mu\delta^2}{h^2}(u^n_{i+\frac{1}{2}})-u^n_{i-\frac{1}{2}}\frac{\mu\delta^2}{h^2}(u^n_{i-\frac{1}{2}})\Big)-\frac{\nu^2\tau}{2}\frac{\delta^4}{h^4}u^n_i+\frac{\nu\tau}{2}\frac{\mu\delta^3}{h^3}\big(\frac{u^2}{2}\big)^n_i=0
\label{eqn:inv2}
\end{eqnarray}
} \normalsize \noindent and the differential approximation can be
written as: \footnotesize{
\begin{eqnarray}
\displaystyle {\mathcal{P}(x,t,u,\nu,u_x,u_t,u_{xx})=u_t+u\
u_x-\nu\ u_{xx}+(C u_x)_x=0} \label{eqn:diff-approx-burgers}
\end{eqnarray}
} \normalsize \noindent The von Neumann stability analysis of
scheme (\ref{eqn:inv2}) under a linearized form provides the
following necessary conditions for $S$, $CFL$ and
$\Omega_\tau=\Omega\tau$: \footnotesize{
\begin{eqnarray}
\displaystyle CFL^2-2S-2\Omega_\tau\leq0,\ \ \ \ 0\leq(4 S)/3-2
S^2+\Omega_\tau\leq1/2 \label{eqn:stab-cond}
\end{eqnarray}
} \normalsize \noindent If $\Omega$ is sufficiently close to zero,
these conditions become then sufficient for the linear
formulation.

\subsection{Numerical application}

\noindent The artificial viscosity term $(C u_x)_x$ is build so as to preserve the symmetries of
the Burgers equation. $C$ is function of $(x,t,u,\tau,h)$.\\

\noindent The Burgers equation is solved numerically for the semi-invariant scheme, and the classical
schemes mentioned above. The numerical solutions
are determined in the frame $(F1)$, with an uniform and orthogonal mesh, in the frame $(F2)$ and $(F3)$,
which are respectively obtained with the Galilean transformation $(x,t,u,\nu) \longmapsto (x+0.25\ t,t,u+0.25,\nu)$
and a higher translation velocity $(x,t,u,\nu) \longmapsto (x+0.5\ t,t,u+0.5,\nu)$.
In $(F2)$ and $(F3)$ the orthogonality of the mesh is not preserved. The parameters, which control the
numerical stability, i.e. the CFL number and the cell Reynolds number, are preserved by the Galilean
transformation. Any change of the frame does not modify their value. These parameters satisfy the
linear stability conditions, which are sufficient for all the numerical applications realized below.\\

\begin{figure}[!h]
\resizebox*{0.49\columnwidth}{0.20\textheight}{\includegraphics{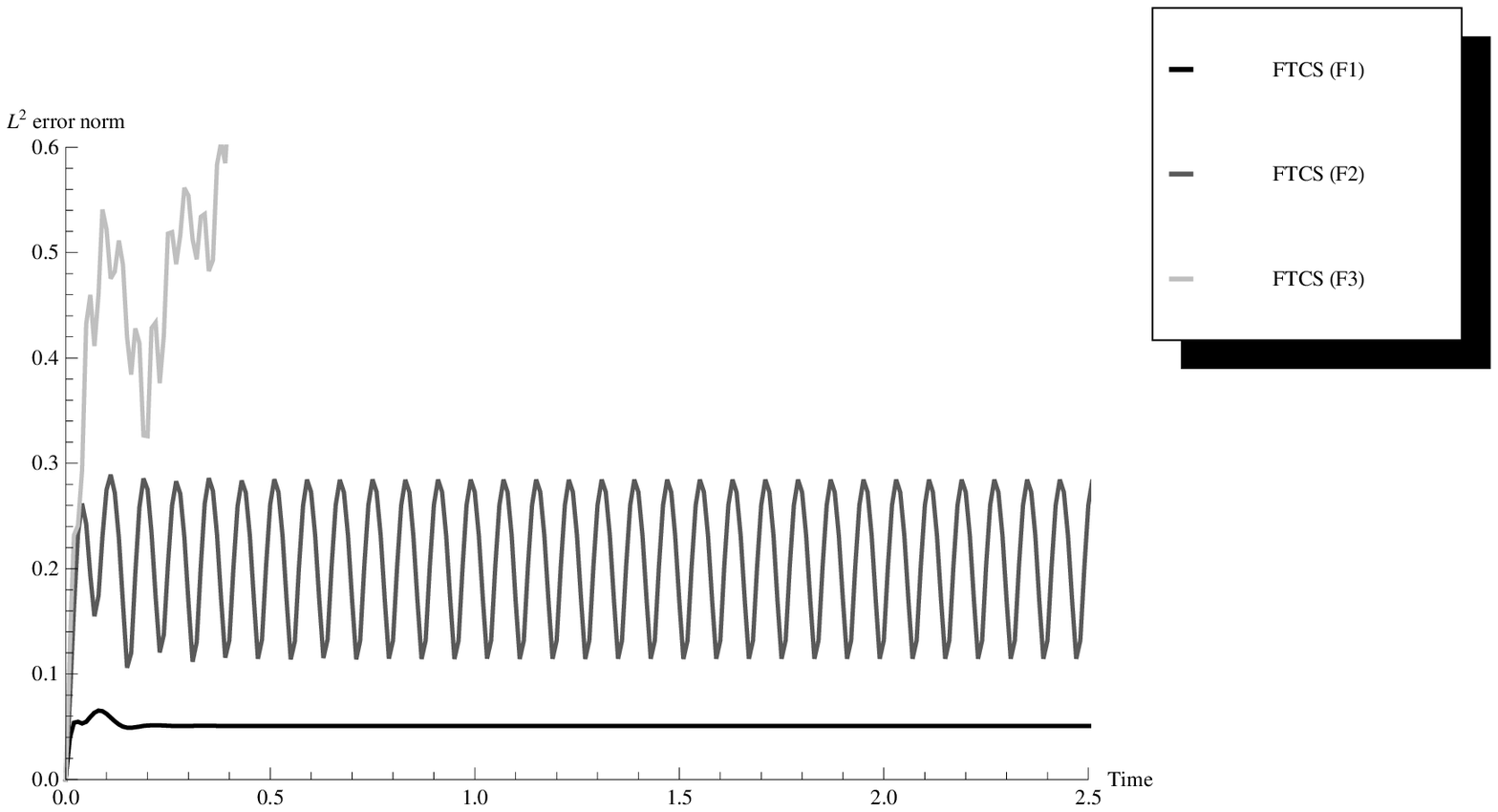}}
\resizebox*{0.49\columnwidth}{0.20\textheight}{\includegraphics{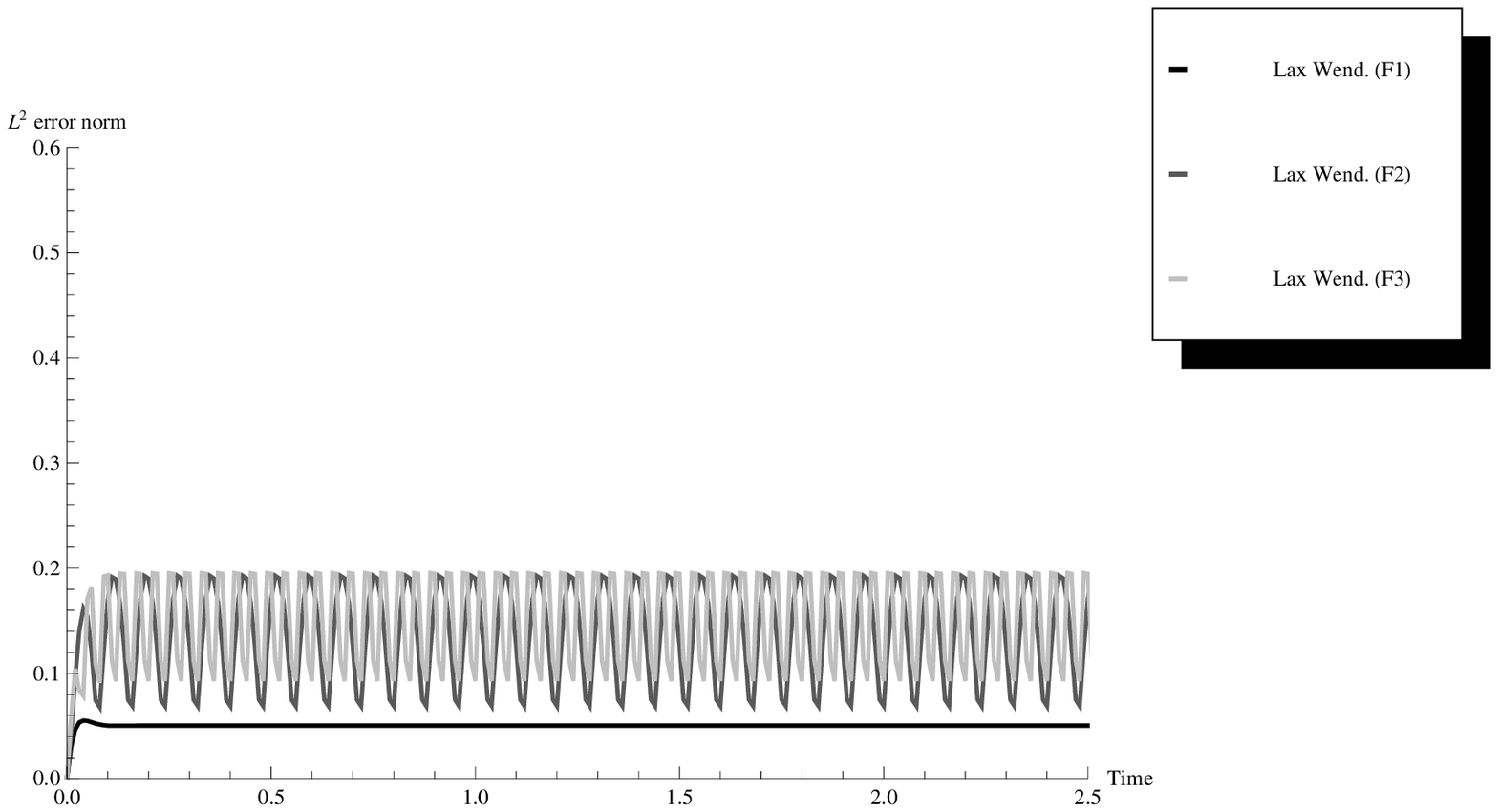}}\\
\resizebox*{0.49\columnwidth}{0.20\textheight}{\includegraphics{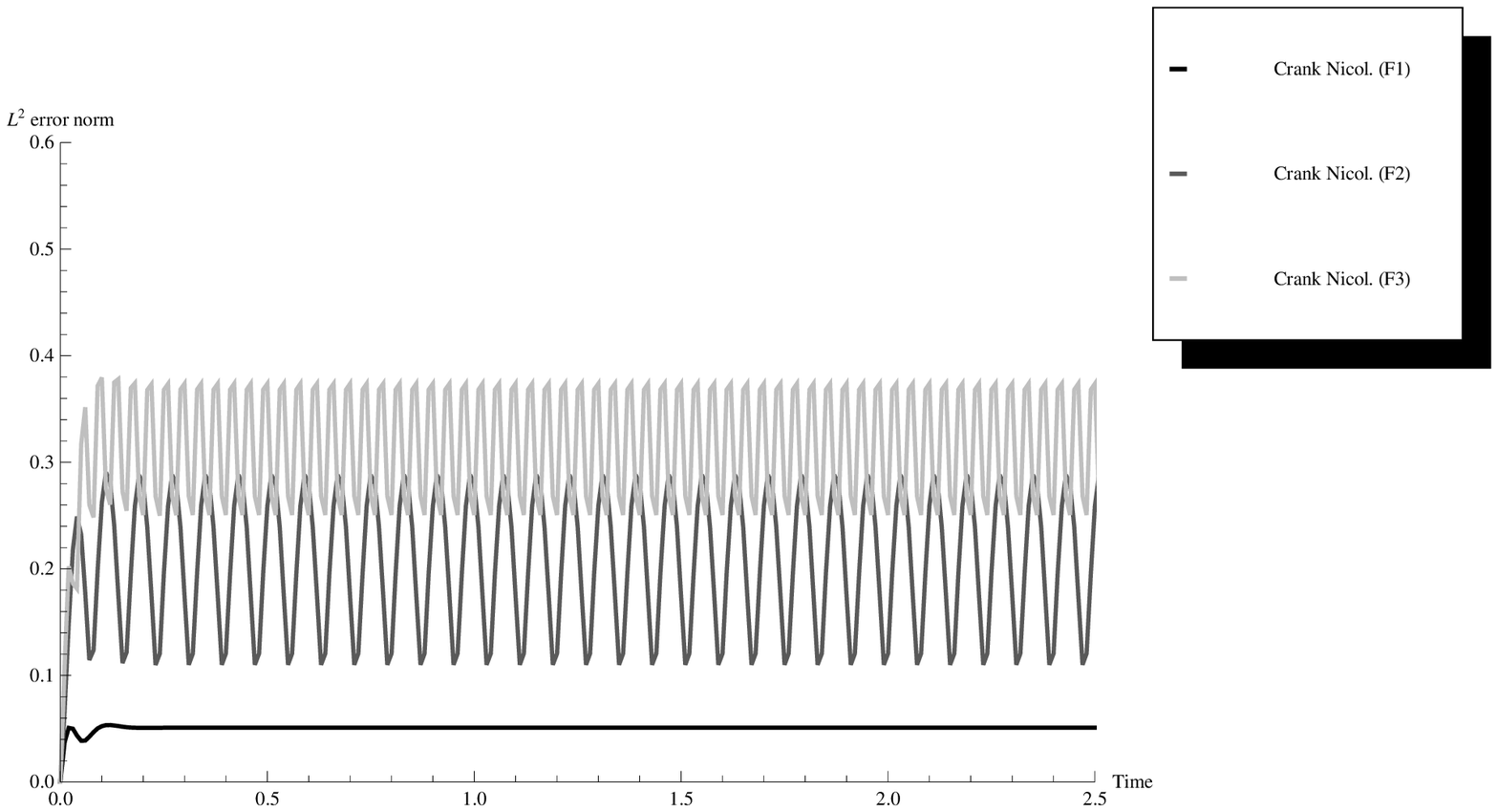}}
\resizebox*{0.49\columnwidth}{0.20\textheight}{\includegraphics{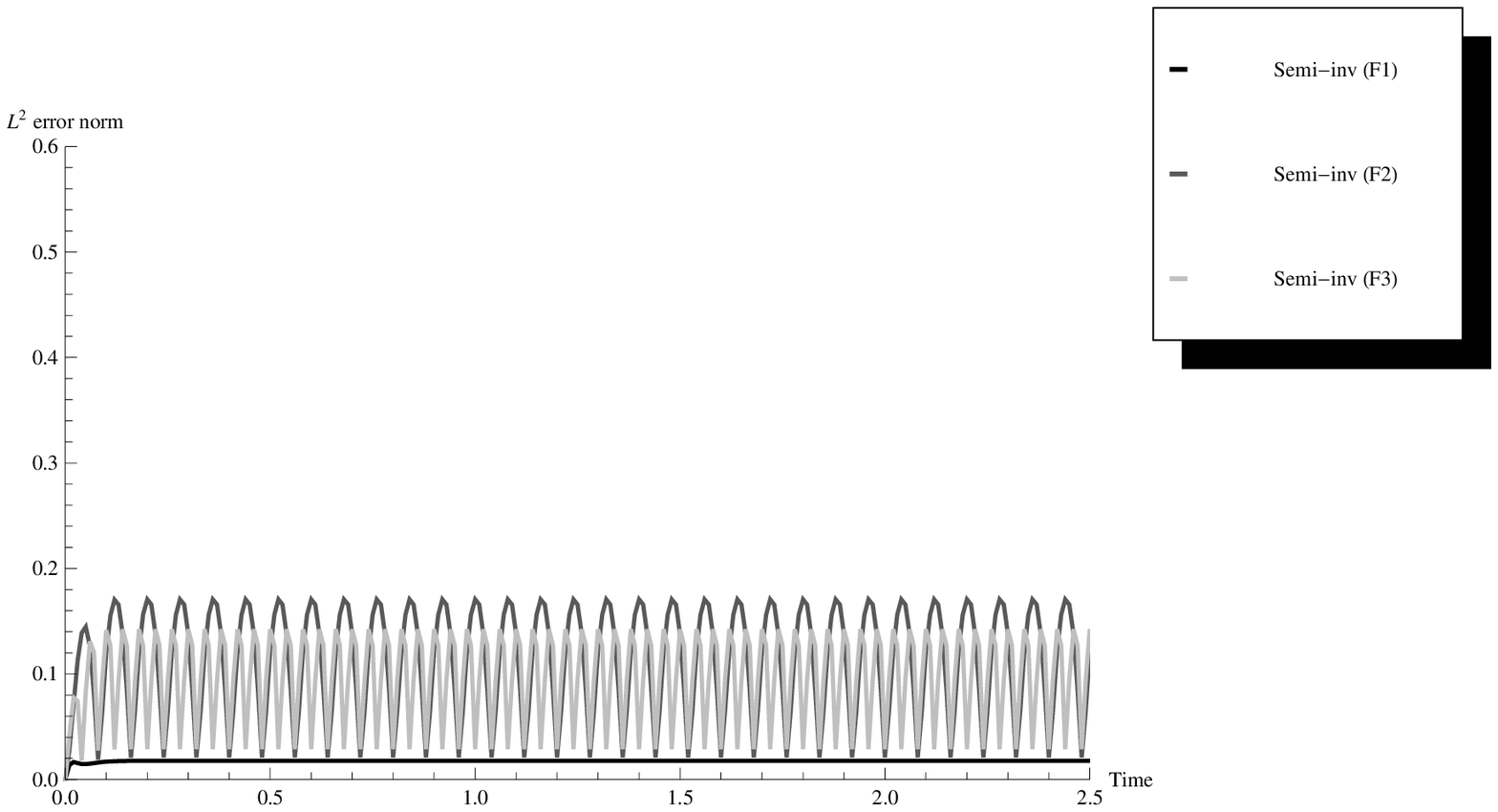}}
\caption{\tiny{Evolution in time of the $L^2$ norm of the error in (F1), (F2), (F3)
when $Re_h=4.0$, $CFL=0.5$} \normalsize} \label{L2-rough-schemes}
\end{figure}

\noindent Figures \ref{L2-rough-schemes} display the time
evolution, of the $L^2$ norm of the error, for the above numerical
schemes, in the different frames. A shock wave type solution is
initially introduced. The CFL and the cell Reynolds number are
respectively equal to $4$ and $0.5$. Non-invariant schemes are
sensitive to the change of frame. Thus, the numerical error of the
FTCS scheme oscillates after the application of a translation with
a velocity of $0.25$. It increases sharply after a translation
with a velocity of $0.5$. The linear numerical stability is
satisfied, though the discretization parameters of the CFL and the
cell Reynolds numbers are situated on the boundary between the
region where the solutions rapidly blow up and the one where the
solutions remain bounded for certain time steps. The introduction
of an additional error due to the loss of the Galilean
transformation leads to the numerical solution blow-up.

\noindent We notice the dependence of numerical error, for the semi-invariant scheme, upon the frame.
The change of the frame causes the apparition of slight oscillations.

\begin{figure}[!h]
\resizebox*{0.49\columnwidth}{0.20\textheight}{\includegraphics{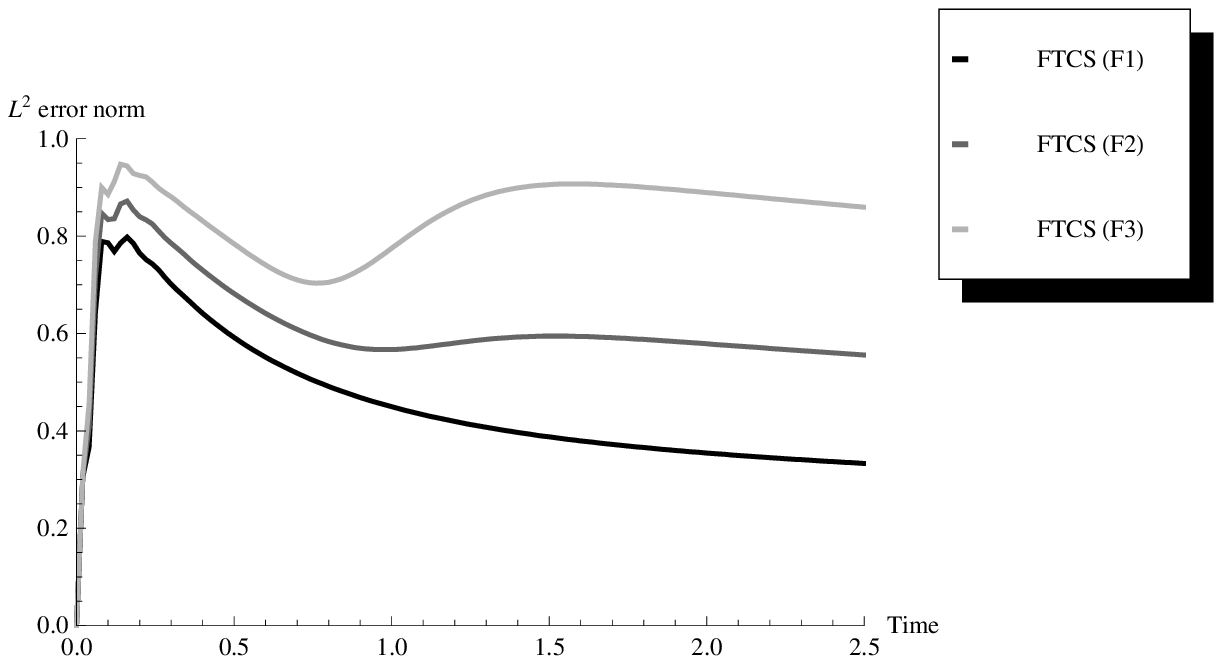}}
\resizebox*{0.49\columnwidth}{0.20\textheight}{\includegraphics{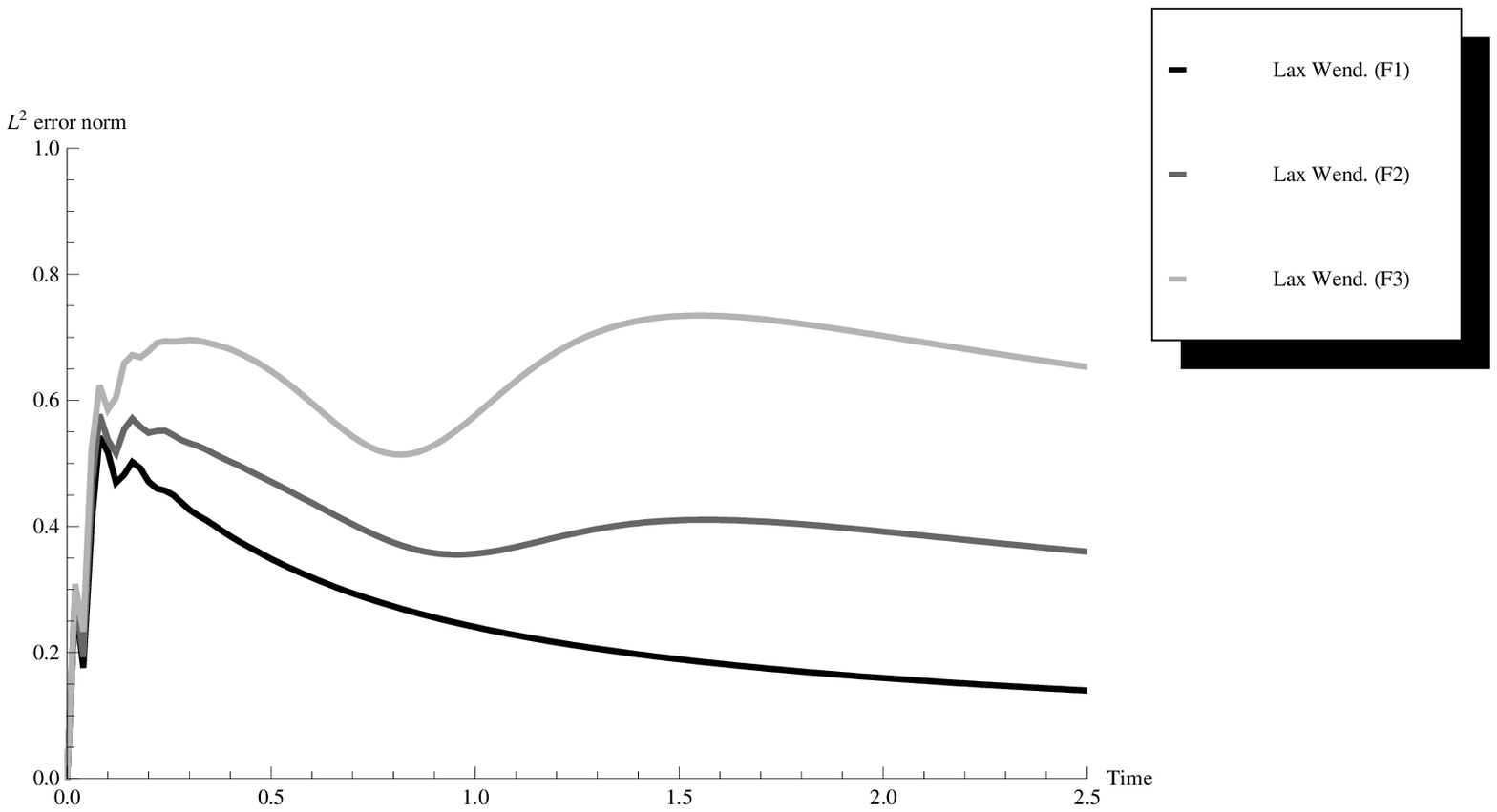}}\\
\resizebox*{0.49\columnwidth}{0.20\textheight}{\includegraphics{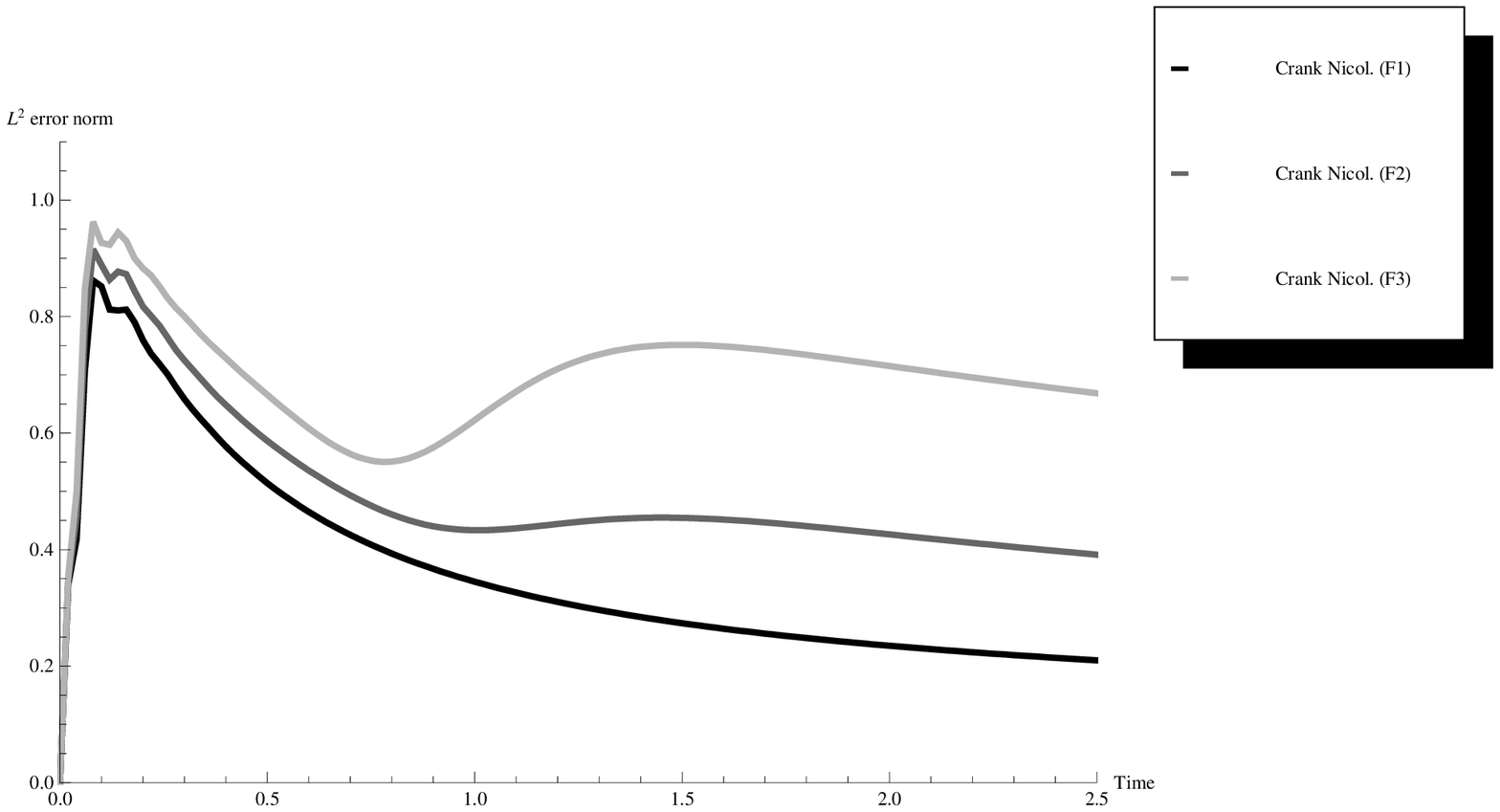}}
\resizebox*{0.49\columnwidth}{0.20\textheight}{\includegraphics{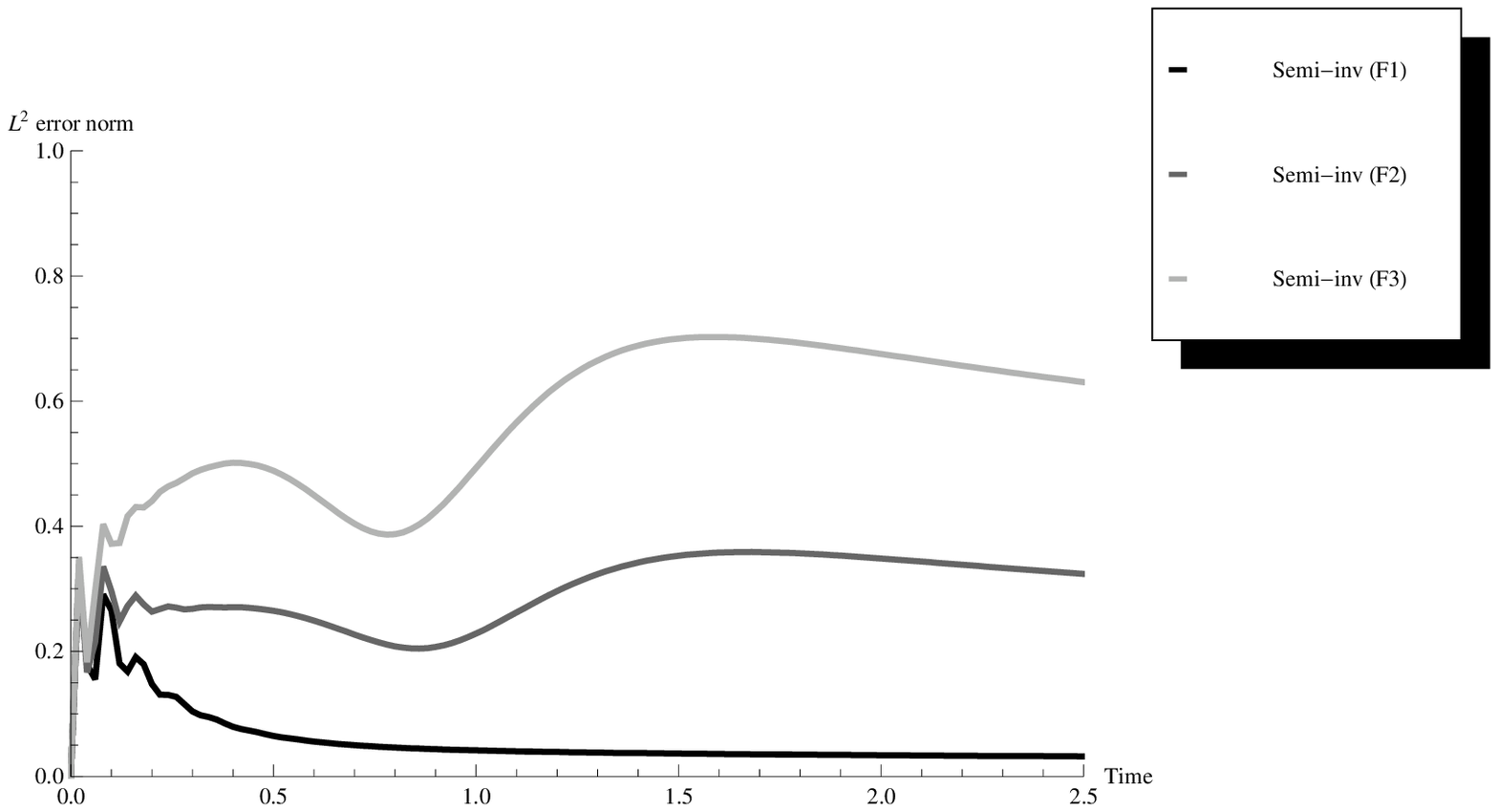}}
\caption{\tiny{Evolution in time of the $L^2$ norm of the error in (F1), (F2), (F3)
when $Re_h=$, $CFL=$} \normalsize} \label{L2-smooth-schemes}
\end{figure}

\noindent Figures \ref{L2-smooth-schemes} display the time evolution, of the $L^2$ norm of the error,
for the above numerical schemes, in the same frames. This time, a wavy solution is initially introduced.
The CFL number and the cell Reynolds number are respectively equal to $0.2$ and $0.06$. We clearly notice
the dependence of numerical error, for the four numerical scheme, upon the frame.

\section{Conclusion}

\noindent Starting from the general \textit{CBKDV} equation, we have set methods that enable us to
build new numerical schemes. The application to the Burgers equation can be seen as a restricted
case of the \textit{CBKDV} equation.

\noindent The invariance method based on the differential
approximation does not enable us to build schemes, which are
rigorously invariant under the action of the group of the original
equations. Like in the case of classical finite difference methods
for aerodynamics, the method based on the differential
approximation breaks the symmetries, which act on the independent
variables, and transform the geometrical properties of the mesh.
This symmetry action is observed for the Galilean transformation
and has been applied for the numerical resolution of the Burgers
equation. Under the action of the Galilean transformation, a
non-invariant discretization leads to the loss of the mesh
orthogonality.

\noindent If, on the one hand, the method based on the
differential approximation would provide the exact invariance if
the higher order terms of the numerical error were considered, it
would yield, on the other hand, a scheme too complicated to be
implemented. In order to preserve rigorously the symmetries, we
need to use direct methods of symmetry analysis of finite
difference schemes, which take into account the transformation of
all the finite difference variables and preserve the geometrical
properties of the mesh.



\addcontentsline{toc}{section}{\numberline{}References}


\begin{thebibliography}{1}

\normalsize{

\bibitem{HoarauDavid}
     {E. Hoarau}, \textbf{Cl. David},
     {Lie group computation of finite difference schemes},
     Dynamics of Continuous, Discrete and Impulsive Systems (Series A), {14}(2007), 180--184.

\bibitem{HoarauDavid2}
     {E. Hoarau}, {Cl. David}, P. Sagaut, T.-H. L\^e,
     {Lie group study of finite difference schemes},
     Discrete and Continuous Dynamic Systems, Supplement 2007, 495--505.

\bibitem{Shokin}
     Yu. I. Shokin,
     {The method of differential approximation},
     Springer-Verlag (1983), Berlin, Heidelberg, New-York, Tokyo.

\bibitem{Yanenko}
     N. N. Yanenko and Yu. I. Shokin,
     {Group classification of difference schemes for a system of one-dimensional equations of gas dynamics},
     Amer. Math. Soc. Transl., {2}, no. 104(1976), 259--265.

\bibitem{Olverappl}
     P. J. Olver,
     {Applications of Lie Groups to Differential Equations},
     Springer-Verlag (1986), New-York.

\bibitem{Ibragimov}
     N. H. Ibragimov,
     {CRC Handbook of Lie Group Analysis of Differential Equations},
     Vol. 1, 2, 3. CRC Press (1994), 1994--1996.

\bibitem{Cantwell}
     B. J. Cantwell,
     {Introduction to symmetry analysis},
     Cambridge University Press (2002).

\bibitem{Ames}
     W. F. Ames, F. V. Postell and E. Adams,
     {Optimal numerical algorithms},
     Applied Numerical Mathematics, {10}(1992), 235--259.

\bibitem{Olverinv}
     P. J. Olver,
     {Geometric foundations of numerical algorithms and symmetry},
     Appl. Alg. Engin. Comp. Commun, {11}(2001), 417--436.

\bibitem{Kim}
     P. Kim,
     {Invariantization of numerical schemes using moving frames},
     Mathematical Physics Seminar (2004).

\bibitem{Dorodnitsyn-theory}
     V. A. Dorodnitsyn,
     {Transformation groups in the net spaces},
     J. Sov. Math., {55}(1991), 1490--1517.

\bibitem{Dorodnitsyn-model}
     V. A. Dorodnitsyn,
     {Finite difference models entirely inheriting continuous symmetry of original differential equations},
     Int. J. Mod. Phys., Series C, {5}, no. 4(1994), 723--734.

\bibitem{Dorodnitsynal2}
     V. A. Dorodnitsyn and R. Kozlov,
     {A heat transfer with a source: the complete set of invariant difference schemes},
     J. Non. Math. Phys., {10}, no. 1(2003), 16--50.

\bibitem{Bakirova}
     M. I. Bakirova, V. A. Dorodnitsyn and R. Kozlov,
     {Symmetry-preserving difference schemes for some heat transfer equations},
     J. Phys. A: Math. Gen., {30}(1997), 8139--8155.

\bibitem{Valiquette}
     F. Valiquette and P. Winternitz,
     {Discretization of partial differential equations preserving their physical symmetries},
     J. Phys. A: Math. Gen., {38}(2005), 9765--9783.

\bibitem{VuCarminati}
     K. Vu and J. Carminati,
     {Symbolic computation and differential equations: Lie symmetries},
     J. Symbolic Computation, {29}(2000), 95--116.

\bibitem{Herod}
     S. Herod (1992),
     {MathSym: a Mathematica program for computing Lie symmetries},
     Preprint, Program in Applied Mathematics, Boulder, Colorado, The University of Colorado.

\bibitem{Baumann}
     G. Baumann (1992),
     {Lie symmetries of differential equations: A mathematica program to determine Lie symmetries},
     Wolfram Research Inc., Champaign, Illinois, MathSource 0202-622.

\bibitem{Schwarz}
     F. Schwarz (1982),
     {A REDUCE package for determining Lie symmetries of ordinary and partial differential equations},
     Comput. Phys. Commun., {27}, 179--186.

\bibitem{Hildebrand}
     F. B. Hildebrand (1956),
     {Introduction to Numerical Analysis},
     New York: McGraw-Hill.

\bibitem{burger1}
     J. M. Burgers,
     {Mathematical examples illustrating relations occurring
     in the theory of turbulent fluid motion},
     Trans. Roy. Neth. Acad. Sci., Amsterdam, {17}(1939), 1--53.

\bibitem{kdv1}
     D. J. Korteweg and G. de Vries,
     {On the change of form of long waves advancing in a
     rectangular channel, and on a new type of long stationary waves},
     Phil. Mag., {39}(1895), 422--443.

\bibitem{wadati}
     M. Wadati,
     {The modified Korteweg-de Vries equation},
     J. Phys. Soc. Japan, {34}(1973), 1289--1296.

\bibitem{wang1}
     M. L. Wang,
     {Exact solutions for a compound KdV-Burgers equation},
     Phys. Lett. A, {213}(1996), 279--287.

\bibitem{feng1}
     Z. Feng and G. Chen,
     {Solitary Wave Solutions of the Compound Burgers-Korteweg-de Vries Equation},
     Physica A, {352}(2005), 419--435.

\bibitem{feng2}
     Z. Feng,
     {A note on ``Explicit exact solutions to the compound Burgers-Korteweg-de
      Vries equation"},
      Phys. Lett. A, {312}(2003), 65--70.

\bibitem{feng3}
     Z. Feng,
     {On explicit exact solutions to the compound Burgers-KdV equation},
     Phys. Lett. A, {293}(2002), 57--66.




\bibitem{David,Feng} {Cl. David}, R. Fernando, Z. Feng,{ A note on "general solitary wave solutions
of the Compound Burgers-Korteweg-de Vries Equation"}, Physica A:
Statistical and Theoretical Physics, {375}(1)(2007), 44--50.



\bibitem{park1}
     E. J. Parkes and B. R. Duffy,
     {Traveling solitary wave solutions to a compound KdV-Burgers equation},
     Phys. Lett. A, {229}(1997), 217--220.

\bibitem{park2}
     E. J. Parkes,
     {A note on solitary-wave solutions to compound KdV-Burgers equations},
     Phys. Lett. A, {317}(2003), 424--428.

\bibitem{zhang1}
     W. G. Zhang, Q. S. Chang and B. G. Jiang,
     {Explicit exact solitary-wave solutions for compound KdV-type and compound
     KdV-Burgers-type equations with nonlinear terms of any order},
     Chaos, Solitons \& Fractals, {13}(2002), 311--319.

\bibitem{zhang2}
     W. G. Zhang,
     {Exact solutions of the Burgers-combined KdV mixed equation},
     Acta Math. Sci., {16} (1996), 241--248.

\bibitem{li}
     B. Li, Y. Chen Y. and H. Q. Zhang,
     { Explicit exact solutions for new general two-dimensional KdV-type
     and two-dimensional KdV-Burgers-type equations with nonlinear terms of any order},
     J. Phys. A (Math. Gen.), {35}(2002), 8253--8265.

\bibitem{whitham}
     G. B. Whitham,
     {Linear and Nonlinear Wave},
     Wiley-Interscience (1974), New York.

\bibitem{ablowitz}
     M. J. Ablowitz and H. Segur,
     {Solitons and the Inverse Scattering Transform},
     SIAM (1981), Philadelphia.


\bibitem{dodd}
     R. K. Dodd, J. C. Eilbeck, J. D. Gibbon and H. C. Morris,
     {Solitons and Nonlinear Wave Equations},
     London Academic Press (1983), London.

\bibitem{johnson}
     R. S. Johnson,
     {A Modern Introduction to the Mathematical Theory of Water Waves},
     Cambridge University Press (1997), Cambridge.

\bibitem{ince}
     E. L. Ince (1956),
     {Ordinary Differential Equations},
     Dover Publications, New York.

\bibitem{zhang3}
     Z. F. Zhang , T. R. Ding, W. Z. Huang and Z. X. Dong,
     {Qualitative Analysis of Nonlinear Differential Equations},
     Science Press (1997), Beijing.

\bibitem{birk}
     G. Birkhoff and G. C. Rota, {Ordinary Differential Equations},
     Wiley (1989), New York.
}

\end{thebibliography}
\end{document}